\documentclass{amsart}
\usepackage{graphics,verbatim,color,amsthm}

\usepackage[all]{xy}

\theoremstyle{plain}
\newtheorem{theorem}{Theorem}

\theoremstyle{remark}

\theoremstyle{definition}

 \def\into{\rightarrow}

\def\R{\mathbb{R}}

\begin{document}

\title[Free homotopy classes for the Three-Body Problem]{Realizing all Free homotopy classes for the Newtonian Three-Body Problem}
\markright{Three-body}
\author{Richard Moeckel}
\address{School of Mathematics\\ University of Minnesota\\
Minneapolis MN 55455}
\email{rick@math.umn.edu} 

\author{Richard Montgomery}
\address{Mathematics Department\\ University of California, Santa Cruz\\
Santa Cruz CA 95064}
\email{rmont@ucsc.edu}

\date{November 30, 2014}

\maketitle

\begin{abstract}  The configuration space of the planar three-body problem when  collisions are excluded 
has  a rich topology  which supports  a large set of free homotopy classes.  Most  classes survive   modding out by
rotations. Those that survive are called the reduced free homotopy classes and have a simple description
when projected onto the shape sphere. They are coded by syzygy sequences.
We prove that every reduced free homotopy class, and thus  every reduced syzygy sequence, is realized by
a reduced periodic solution to the Newtonian planar three-body problem.  The realizing solutions
have nonzero angular momentum,  repeatedly come very close to triple collision, and have lots of ``stutters"--
repeated syzygies of the same type.  The heart of the proof is contained in the   work by one of us    on  
symbolic dynamics arising out of the central configurations after the triple collision is blown up using
McGehee's method. 
\end{abstract}




\section{Introduction}
A basic theorem in Riemannian geometry inspires our work. Recall   two loops in a space
$M$ are   {\it freely homotopic} if one loop can be deformed into the other without leaving $M$. 
The resulting equivalence classes of loops are the {\it free homotopy classes}.   This basic theorem asserts
if $M$ is a compact
Riemannian manifold then   every one of its free homotopy classes of loops is realized by   a   periodic geodesic. 

This theorem  suggests an analogue for the  planar  Newtonian three-body problem.   Replace the Riemannian manifold  above  by    the configuration space  $M$ of the planar three-body problem: the product of $3$ copies of the plane, minus collisions.    Is every free homotopy class
of  this $M$  realized by a solution to the planar three-body problem?  By a {\it reduced free homotopy class} we mean
a free homotopy class of loops  for the  quotient space $M/SO(2)$ of  the configuration space $M$ by the  group $SO(2)$ of rotations  acting on
$M$  by  rigidly  rotating the   triangle formed by the three bodies.   By a {\it reduced periodic solution} we mean a solution which is periodic modulo
rotations, or, what is the same thing, a solution which is periodic in some rotating frame. 
 \begin{theorem}
 \label{th_1}
Consider the planar three-body problem with fixed negative energy and either equal or near equal masses.  Then, for all sufficiently small {\em nonzero} angular momenta,
 every reduced free homotopy class is  realized by a reduced periodic solution. 
\end{theorem}
The case of zero or large angular momentum remains open.

For further motivation and history regarding this problem please see \cite{Albouy_Cabral}. 


{\bf Syzygy Sequences.} 
The reduced free homotopy class of a periodic  motion of the three bodies  can be read off  of its  {\it syzygy sequence}. 
A `syzygy' is a collinear configuration   occurring  as
the three bodies move. (The word comes from  astronomy.)  
Syzygies  come in three types, type   1,2, and 3 , depending on the label of the  mass  in the middle at collinearity. 
A  generic curve in $M$ has a discrete set of syzygies.  List its   syzygyies  in temporal order to obtain the  syzygy sequence of that curve.  
If two or more of the same letters occur in a row, for example 11,   cancel them in pairs using  a homotopy. 
We call such multiple sequential occurrences of the same letter a ``stutter".  Continue canceling stutter pairs until there are no more.   For example $121123 \mapsto 1223 \mapsto 13$.  
We call the final result of this cancellation process  the {\em reduced syzygy sequence} of the curve.  
Periodic reduced syzygy sequences are in a 1:2 correspondence with
reduced free homotopy classes.  Theorem  1 implies that
all reduced syzygy sequences are realized by reduced periodic solutions.
 
\begin{figure}[h]
\scalebox{0.6}{\includegraphics{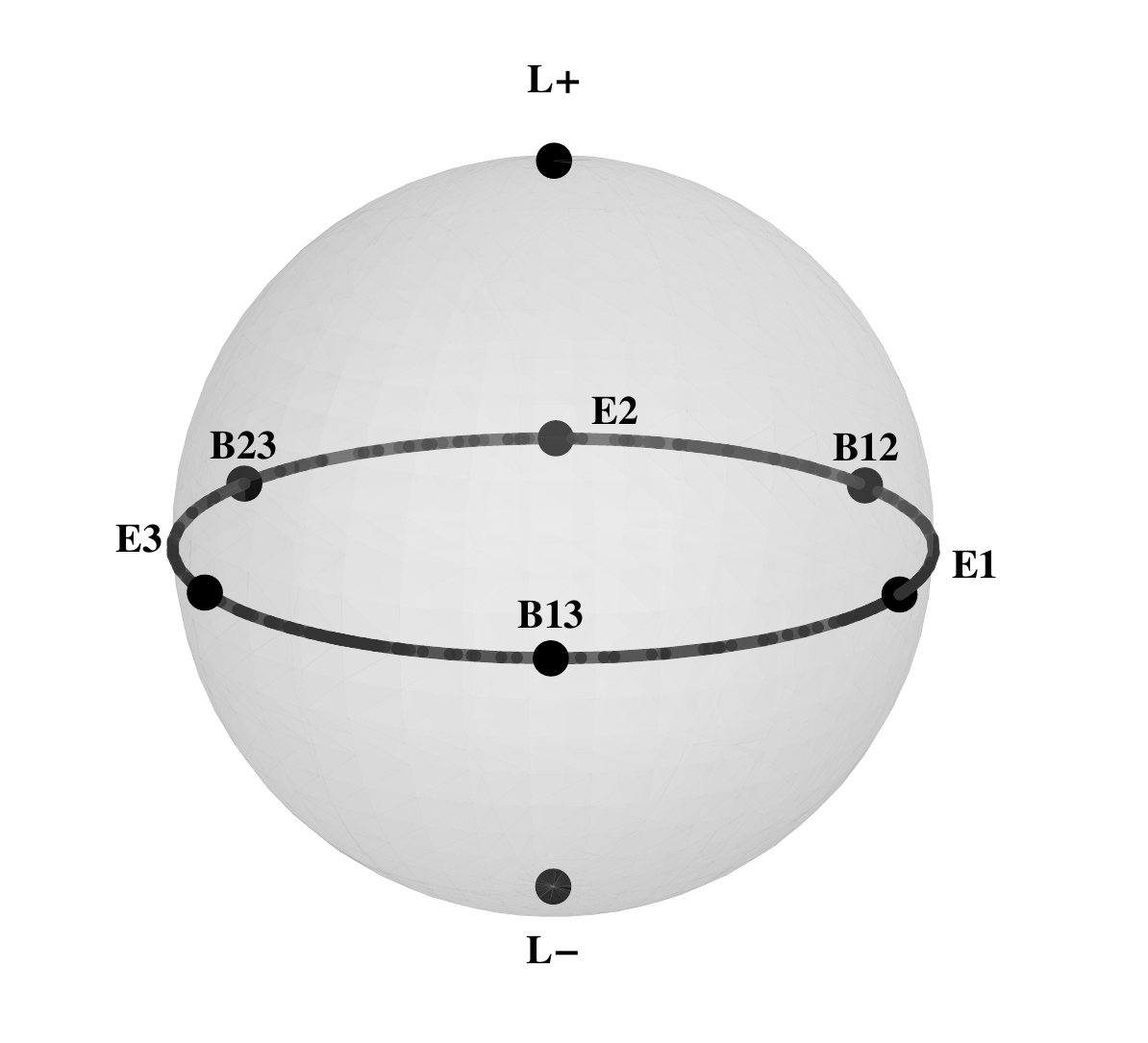}}
\caption{The shape sphere.  The equator represents collinear shapes, including the binary collisions $B_{ij}$ and the Eulerian central configurations $E_i$.  The poles represent the Lagrangian, equilateral central configurations $L_\pm$.}  
\label{fig_shapesphere1}
\end{figure}

 {\bf Shape space and topology.} 
After fixing the center of mass, the configuration space $M$ of the planar three-body problem  
is diffeomorphic to the product $\R^+\times S^1 \times (S^2 \setminus \{ \text{ three points } \})$ where the first factor represents the size of the triangle, the second an overall rotation and the third the shape of the triangle.   The  rotation group $SO(2)$ generates the circle factor $S^1$ and can be identified with it.
 We will forget the circle of rotations,  and focus on homotopy classes and solutions {\it modulo} rotation.   Up to homotopy equivalence, the size factor is also irrelevant.
 The  remaining  $S^2 \setminus \{ \text{ three points } \}$  is
 the {\it shape sphere} minus its three  binary collision points.   See figure \ref{fig_shapesphere1}.   Points of the shape sphere represent
 oriented similarity classes of triangles.  There is a canonical quotient map   $M \to S^2 \setminus \{ \text{ three points } \} $
 which  sends a configuration to its {\it shape}. 
 The equator  of the shape sphere represents   collinear configurations: all three masses in a line.
The  three deleted collision points, denoted $B_{12}, B_{23}, B_{13}$  lie on this  equator and  split it into three arcs, which could be labelled 1, 2, and 3 according to the syzygy type.

\begin{figure}
\scalebox{0.6}{\includegraphics{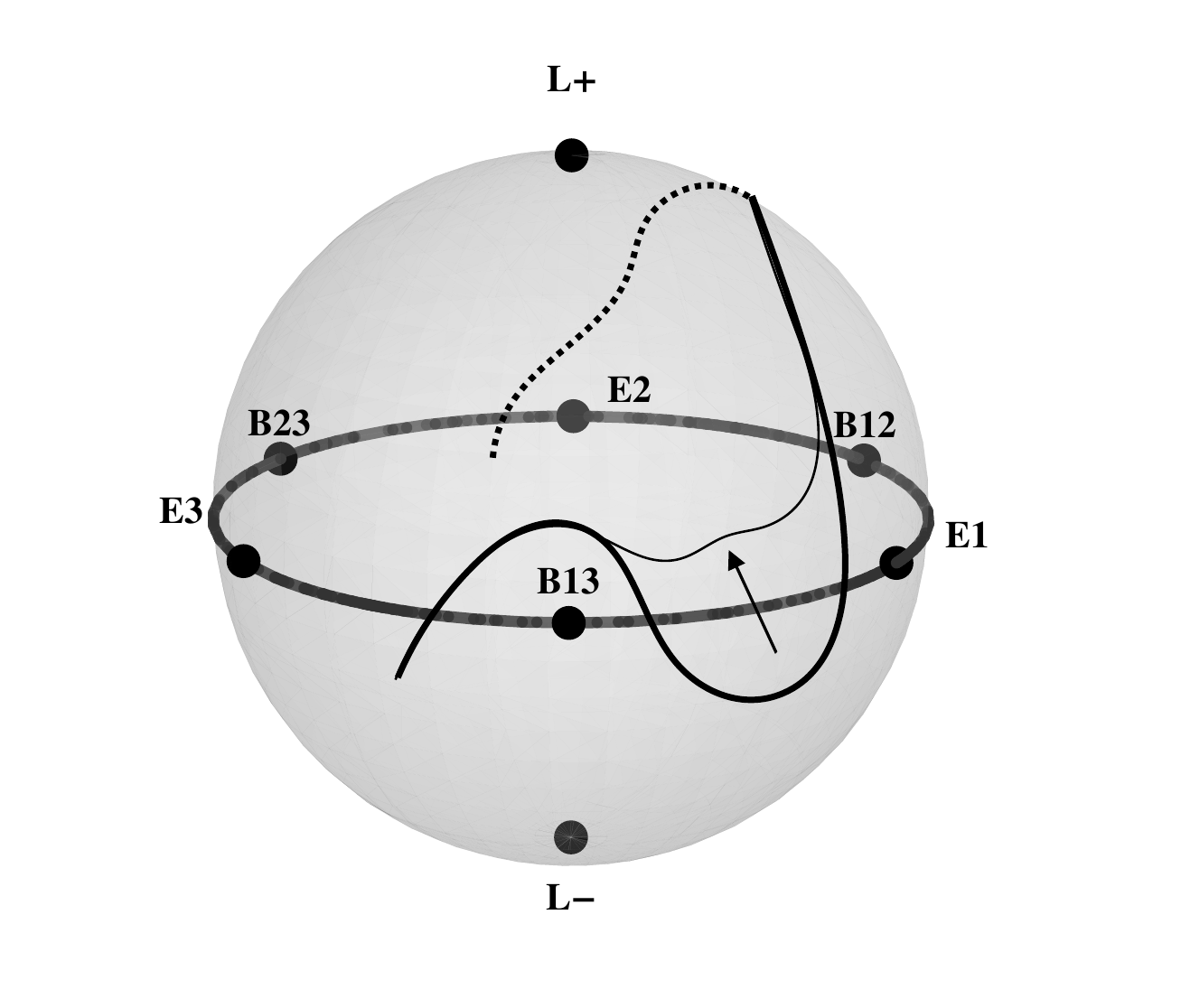}}
\caption{Canceling Stutters.  A path with syzygy sequence 3112 is homotopic with fixed endpoints to one with sequence 32.}  
\label{fig_stutter}
\end{figure}

We redo the construction of the  syzygy sequence  of a curve  using shape language.
Take a  closed curve in $M$. Project it to the shape sphere, perturbing it if necessary so that it crosses the equator transversally.   Make a list of the arcs encountered in temporal order.  This list is the curve's syzygy sequence.  
Refer to figure \ref{fig_stutter}  for the picture in the shape sphere of cancelling stutter pairs.
We cancel stutter pairs until there are no more, arriving at the  reduced syzygy sequence of the loop. 
The   reduced syzygy sequence  is a free homotopy invariant:   two loops, freely homotopic in $M/SO(2)$,   have the same reduced syzygy sequence.

This map from reduced free homotopy classes to  reduced syzygy sequences is two-to-one, with a single exception.
The  two  reduced free homotopy classes which represent a given  nonempty reduced syzygy sequence are
related  by reflection, or, what is the same thing, by the choice of direction ``north-to-south" or ``south-to-north"  with
which we   cross the equator when we list the first letter of our syzygy  sequence. 
The single exception to this two-to-one rule  is  the empty 
reduced syzygy sequence having no letters.  This empty class represents only the trivial reduced  free homotopy class whose
representatives are the contractible loops in $M/SO(2)$. 

Theorem \ref{th_1} is an immediate corollary of the following theorem.
To state it,  define a stutter block of size $n$ to be a syzygy sequence of the form $\epsilon^n$ where $\epsilon\in\{1,2,3\}$.
The theorem will guarantee that we can realize any syzygy sequence which is a concatenation of such stutter blocks with
sizes in a certain range $N\le n\le N'$.    Let's call such a syzygy sequence a {\em stutter sequence with range $[N,N']$}.  Note that one could put several blocks with the same symbol together so while $N$ is a lower bound on the number of repetitions of a symbol, $N'$ is not an upper bound.

\begin{figure}[h]
\scalebox{0.6}{\includegraphics{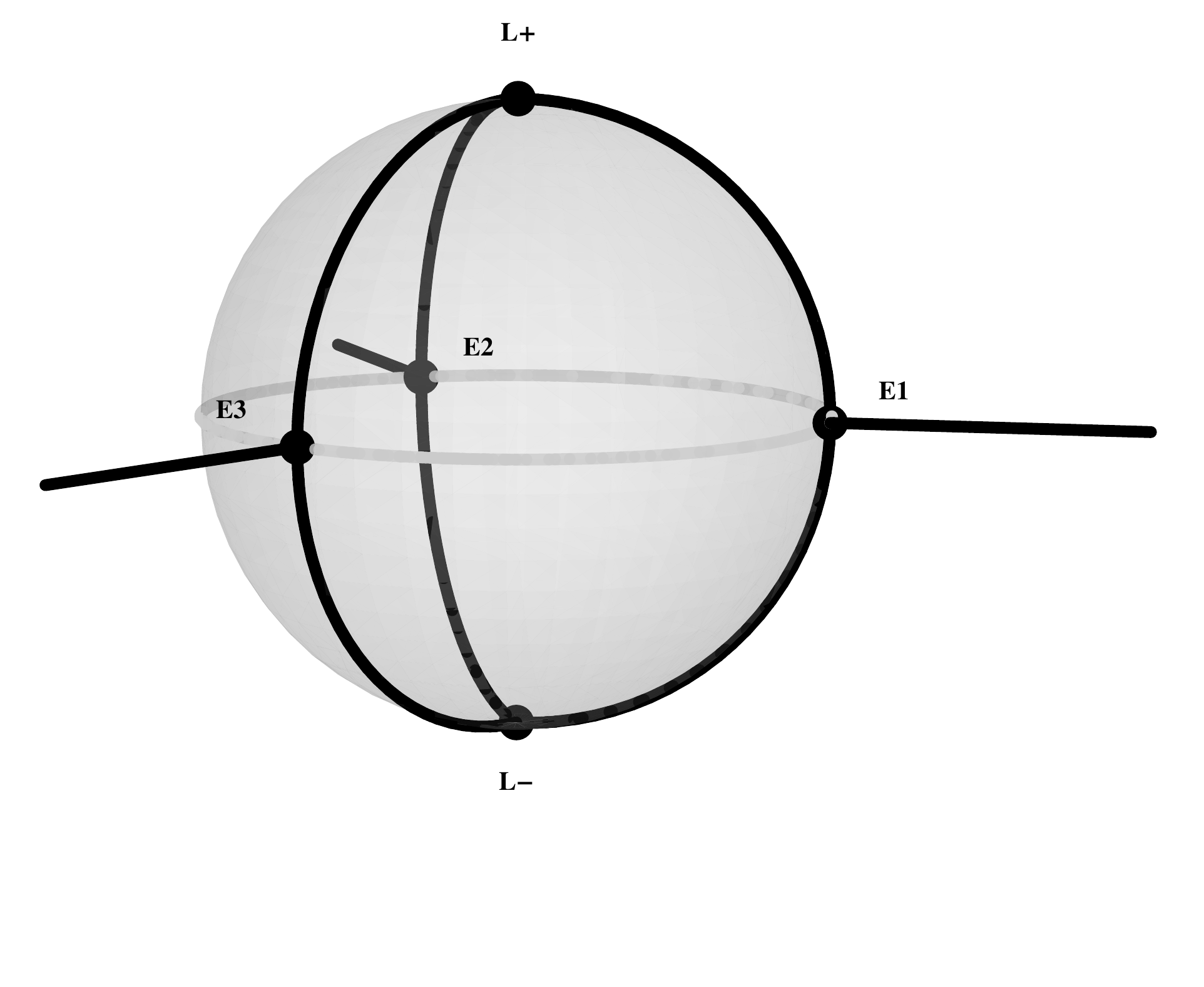}}
\caption{The solutions in the reduced configuration space move near the black curves. The sphere represents configurations with $\{r=0\}$ (triple collision) and the exterior represents $\{r>0\}$. The radial line segments are the three Eulerian homothetic solutions and the circular arcs  on the sphere are isosceles configurations }  
\label{fig_skeleton}
\end{figure}

\begin{theorem}\label{th_2}
Consider the planar three-body problem at fixed negative energy.  
There is a positive integer $N$ such that given any $N'\ge N$, all bi-infinite stutter sequences with range $[N,N']$ are realized by a collision-free solution,
 provided the masses are sufficiently close to equal and the angular momentum $\mu$ satisfies  $0<|\mu|<\mu_0(N')$.
If the syzygy  sequence is periodic then
it is realized by at least one reduced periodic solution.   Between each stutter block the realizing solutions pass near triple collision.  
Their projections to the reduced configuration space lie near the black curves of figure \ref{fig_skeleton}.
\end{theorem}

{\sc Proof of theorem 1 from Theorem 2.} Take $N' > N$ in the theorem \ref{th_2} so that we are guaranteed  that 
there is an odd integer $k$ in the interval $[N, N']$.  Given a nontrivial reduced periodic syzygy sequence, replace
each occurrence of  letter $i$ by $i^k$ to get a   corresponding periodic  unreduced syzygy sequence  satisyfing the conditions of   theorem 2.
Since $k$ is  odd, this sequence reduced to the given reduced periodic sequence. 
 The periodic solution guaranteed by theorem \ref{th_2} then realizes  one of the two free homotopy classes having
the given reduced syzygy sequence.  To realize  the other class    apply a   reflection to this solution.  
 Finally, to   realize the empty  sequence we can simply take the Lagrange solution.  Or we could take all the exponents $j_i$ even: for example all $1^k 1^k$'s,
or repeated $1^{2k} 2^{2k}$.  QED
\vskip .3cm

{\bf A History of  Failed Variational Attempts.} 
The proof of the realization of all free homotopy classes  in the Riemannian case uses    the direct method
of the calculus of variations.  Fix a free homotopy class. Minimize the length over all loops realizing
that class.   Since the manifold is assumed  compact, the minimum is realized and is a geodesic. 

In mechanics the extremals of   the  {\it action} 
yield  solutions to Newton's equations.
It is natural to  try to proceed in the same way. Fix a   reduced syzygy sequence. Minimize the action
over all representatives of that class.  The configuration space is noncompact which complicates the approach.
There are two sources of noncompactness: spatial infinity $r \to \infty$ and collision $r \to 0$,
where $r$ denotes the distance between any two of the three bodies.  The first type
is fairly easy to overcome.  Noncompactness due to the collision is the essential difficulty in applying the direct method.
For  the standard $1/r$-type potential of Newtonian gravity,   there are paths {\it with collision} 
which have finite action.  Through these paths we find, by explicit examples, that while   trying to minimize the action we  may leave the free
homotopy class we started in, through a `pinching off' via collision, in which we leave $M$ and we lessen the action, entering  a different free homotopy class.  
In other words: the infimum of the action over a fixed free homotopy class is typically not achieved, and minimizing sequences tend to collision.
See \cite{Mont_Action} for an assertion that this is probably a ubiquitous problem, suffered when we choose almost any free homotopy class. 

The fix to this collision problem, going back at least to Poincar\`e \cite{Poincare}, is to cheat.  Change   the Newtonian potential to one which blows up like   $1/r^2$ (or even stronger)   as
we approach collisions at $r=0$.  
If the negative of the potential for the force is greater than $C/r^2$ whenever the distance $r$ between two of the bodies
is sufficiently small, then  the action of any path suffering a collision
is infinite.   Such potentials  are called ``strong force". 
Under the strong force assumption each free homotopy type is separated from every other by  an infinite ``action wall".
The direct method works like a charm. Poincare \cite{Poincare}  proved that almost every nontrivial homology class in $M$ is represented by a periodic orbit.
About a century later one of us proved \cite{Mont_strongforce} that  almost every free homotopy class modulo rotations, i.e., that almost every reduced syzygy sequence,  is realized
by a   periodic solution to the planar strong force three-body problem.  (The only ones that are not realized are
those that only involve 2 letters, so : 1212.. , 1313... and 2323... . These spiral in to binary collision between the two involved bodies.)  

 For  years one of us had tried to  realize all  reduced  free homotopy classes  for the honest Newtonian $1/r$ potential by  variational methods.
These efforts  led to interesting discoveries \cite{Chenc_Mont}, \cite{Mont_InfinitelyMany}, \cite{Mont_pants}  but   did not seem to   bring us any  closer to
 solving the original realization problem.   (For more on the various attempts and theorems along the way also see \cite{Albouy_Cabral}.) 
  So it was with great joy  that
 we realized that  ideas developed  by the other author  in the 1980's  \cite{M_chaotic},  summarized and   extended as  Theorem 2,
  would  prove of theorem 1.


{\bf The main idea behind the proof of theorem~\ref{th_2}. }  

The solutions described by  Theorem 2 arise out of chaotic
perturbations of the  homothetic solutions of Euler and Lagrange.   We   recall  that the 
solutions of Euler and Lagrange are   precisely 
those solutions to Newton's equations whose shape does not change as they evolve. When projected onto   the shape sphere such a solution curve is a single constant point.
There are five points in all.  They are called central configurations and denoted here $E_1, E_2, E_3, L_+, L_-$.
The $E_i$ are Euler's collinear configurations.  $E_i$  lies on the collinear arc marked $i$. The $L_{\pm}$ are Lagrange's central configurations and are
equilateral triangles, with one for each possible orientation of a labelled equilateral triangle in the plane.

\begin{figure}[h]
\scalebox{1}{\includegraphics{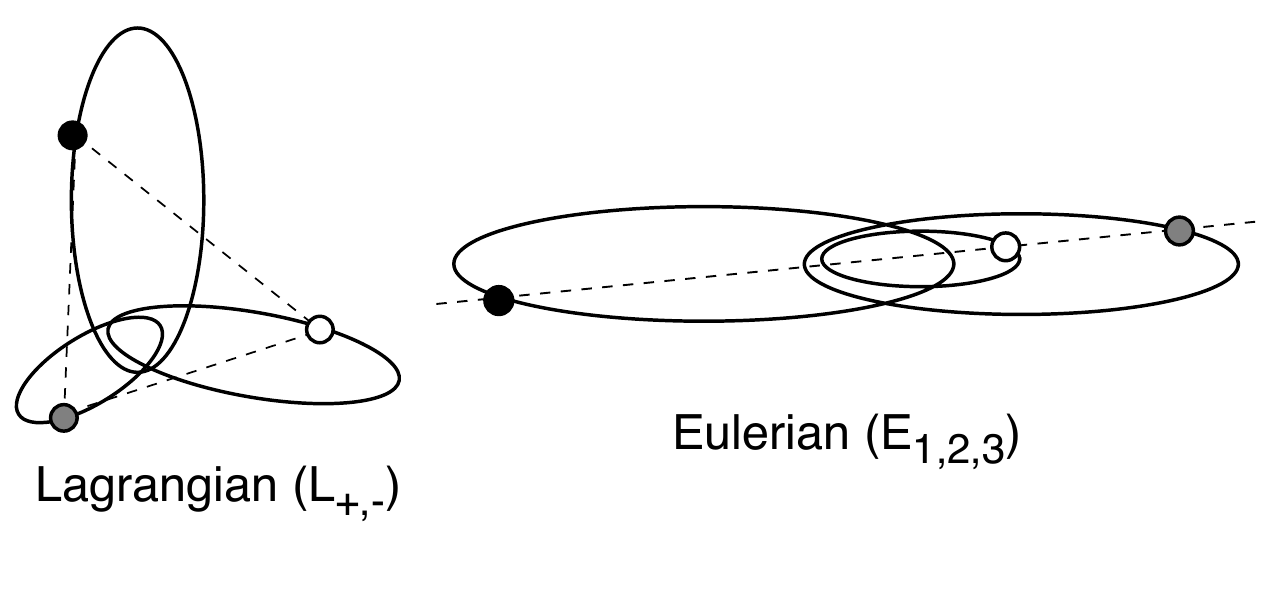}}
\caption{Eccentric homographic solutions}  
\label{fig_eccentriccentral}
\end{figure}

 Each central configuration gives rise to a family of solutions,
the family being parameterized by  an auxiliary Kepler problem, i.e., by  conics.  For each choice of conic the corresponding solution 
consists of the  three masses travelling along  homographic ( =scaled, rotated)   versions of this conic, with the conics  placed so that  one focus is 
 the  center of mass
of the three bodies. The masses travel their
individual conics  in such a way that  the triangle  they form remains in a constant shape.  At fixed energy we can think of the parameter as the angular momentum.
  The homothetic solutions  correspond to angular momentum zero and are degenerate conics. For these solutions the   three bodies
  move along collinear rays, exploding out of triple collision until they reach a maximum size (dependent on the energy) at which instant
  they are all instantaneously at rest.  From that instant they reverse their paths, shrinking back homothetically to triple collision.
  When the angular momentum is turned on to be small but near zero, the masses move along highly eccentric nearly collinear ellipses.
  See figure \ref{fig_eccentriccentral}. 
  
  The solutions of Theorem 2 feature many close approaches to triple collision.  While they are away from collision they look much like the Eulerian homothetic solutions -- the bodies are close to one of the three Eulerian central configurations $E_i$ while the size expands and  then contracts again to another close approach to triple collision.  Near collision, however, their behavior is quite different.  The shape begins to oscillate around the collinear shape, producing a large syzygy block $i^n$.  Very close to collision the shape approaches one of the two Lagrange equilateral triangles $L_\pm$.  The transitions from  $E_i$ to $L_+$ or $L_-$ are through isosceles shapes (the circular arcs in figure~\ref{fig_skeleton}).  One can specify a transition to either $L_+$ or to $L_-$ at each close approach.   After the shape is nearly equilateral, the triangle spins around  (similar to the behavior near collision of the Lagragian homographic solution in figure~\ref{fig_eccentriccentral}).   Next, the shape makes a transition back to one of the  Eulerian central configuration $E_j$, following one of the three isosceles circles.  Arriving near $E_j$ the shape oscillates and the size starts to increase again.  Away from collision, the behavior is like the Euler homothetic solution with shape $E_j$.  The main point is that, during each approach to collision, we can choose which Lagrange configuration to approach and then which Eulerian homothetic orbit to mimic next.  In this way we can concatenate the three types of syzygy blocks at will.

Note that the arrangement of  the three isosceles half-circles in figure~\ref{fig_skeleton}  has the same free homotopy type as the sphere minus three points.
By following the  dark lines of  figure \ref{fig_skeleton} we can realize every  reduced free homotopy class. 
The results of  \cite{M_chaotic} allow us to make the necessary  transitions near the poles from one kind of Eulerian behavior to the other.
What was not explicitly said there, but follows from the proof, is that when the masses are nearly equal these realizing solutions,
when projected to the shape sphere, stay within a small neighborhood of the isosceles half-circles.  This prevents the occurrence of additional, unintended syzygies and also allows us to avoid binary collisions.


There are a number of  reasons we need   some angular momentum to achieve the  gluing of the three half-circles.
Without angular momentum, the isosceles problems are invariant sub-problems and so  we cannot switch from one circle to another.
 Sundman's   theorem   asserts that triple collision is impossible when the angular momentum is not zero,
 hence some   angular momentum relieves us  of the inconvenience of investigating whether or not solutions ``die" in    triple collision. 
But the most important reason is  that we need the Lagrange-like spinning behavior near triple collision to connect up the orbit segments approaching collision with those receding from it.  See figure \ref{fig_loop} where this spinning behavior appears as a restpoint connection at $r=0$.
This mechanism is described in more detail in the next section.

%


\section{Proof of Theorem~\ref{th_2}}  
In this section, we will describe how theorem~\ref{th_2} follows from previous results on constructing chaotic invariant sets for the planar three-body problem with small, nonzero angular momentum \cite{M_chaotic}.  These results, in turn, are based on a study of the isosceles three-body problem, so we begin  there.  

\subsection{Triple collision orbits in the isosceles problem}
\label{subsec_isosc}
When two of the three masses are equal, there is an invariant subsystem of the planar three-body problem consisting of  solutions for which the shape of the triangle formed by the three bodies remains isosceles for all time.  All of these solutions have zero angular momentum.   After fixing the center of mass, the isosceles subsystem has two degrees of freedom so fixing the energy gives a three-dimensional flow.  In the early eighties, there were several papers published about this interesting problem \cite{Devaney, Simo, M_InfinitelyClose}.   We will describe a few features relevant to syzygy problem, referring to \cite{M_InfinitelyClose} for more details.    In order to be definite, let us suppose that we are thinking of $m_1 = m_2$
so that mass 3 forms the vertex of the isosceles triangles.  

The orbits of interest pass near triple collision.  Using McGehee's blow-up method,  we replace the triple collision singularity by an invariant manifold forming a two-dimensional boundary to each three-dimensional energy surface.  The blow-up method involves introducing size and shape variables $r$ and $\theta$ and corresponding momenta $\nu$ and $w$.  The size $r$ of the triangle is the square root of the moment of inertia
with respect to the center of mass.  Thus $r=0$ represents triple collision of the three masses at their center of mass.  A change of timescale by a factor of $r^\frac32$ slows down the orbits near collision and gives rise to a well-defined limiting flow at $r=0$.

The blown-up equations read
\begin{gather*} 
r^{\prime} = r \nu \\ \theta^{\prime} = w \\ \nu^{\prime} = \frac{1}{2} \nu^2 + w^2 - U(\theta) \\ w^{\prime} = \frac{dU}{d\theta} - \frac{1}{2} \nu w 
\end{gather*}
together with the energy constraint
\[\frac{1}{2} \nu^2 + \frac{1}{2} w^2 - U(\theta) = rh \]
Here $\theta \in [-\pi/2, \pi/2]$ parameterizes the isosceles great circle in shape space for the two equal masses, and so describes the shape of
the isosceles triangle.  
The boundary of the $\theta$-interval, $\theta = \pm\frac{\pi}{2}$, represent  binary collision between the two equal masses.   We refer to \cite{Devaney}
for details regarding $U$. 
As $\theta$ runs over $[-\frac{\pi}{2},\frac{\pi}{2}] $, the triangle opens up, becoming equilateral at a certain angle $\theta_-$ and collinear   at $\theta = 0$, then passing through another equilateral shape at $\theta_+$ before collapsing to binary collision again.  
The equilateral and collinear shapes are the Lagrangian ($L_\pm$)  and Eulerian $E_3$ central configurations.

\begin{figure}[h]
\scalebox{0.6}{\includegraphics{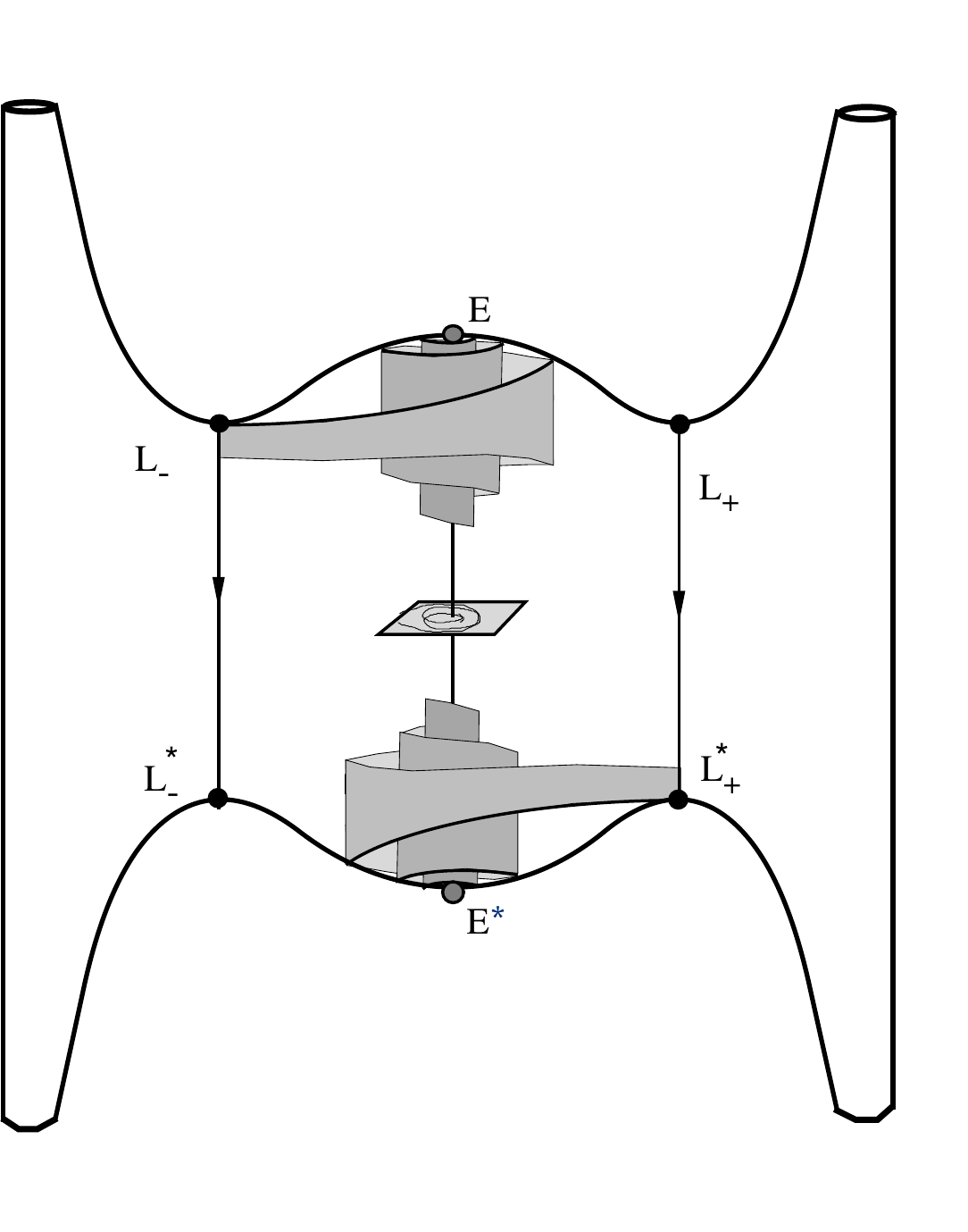}}
\caption{The isosceles three-body problem.  The coordinates are $\theta$ (left to right), $\nu$ (bottom to top) and  $w$ (back to front)}  
\label{fig_isosceles}
\end{figure}

The  triple collision orbits of the three-body problem, isosceles or not,   converge  to a central configuration.  
For each central configuration, there are two restpoints on the associated full collision manifold, one associated to solutions converging to that limiting shrinking shape, the other exploding out of it.
 The isosceles
collision manifold forms  an invariant submanifold of the full planar collision manifold,
 discussed in the next subsection.  For the isosceles subproblem we   have
a total of six rest points  denoted $L_+ ^*, L_- ^*$,  $E_3 ^*$ and $L_+, L_-$ and $E_3$.    (The notation is different in \cite{M_InfinitelyClose}.)
 
 Isosceles solutions having triple collision in forward time   converge  to one of the restpoints $E_3^*,L_-^*,L_+^*$ as the rescaled time tends to infinity.  In other words, they form the stable manifolds of these three restpoints.   Similarly, orbits colliding in backward time (ejection orbits) are represented by the unstable manifolds of $E_3,L_-,L_+$.  In the figure, there are connecting orbits running from the ejection restpoints to the corresponding collision restpoints.  These are the well-known homothetic solutions where the shape remains constantly equal to one of the central configurations while the size expands from $r=0$ to some maximum and then contracts to zero again.  
 
 The six restpoints in the collision manifold are all hyperbolic and the dimensions of their stable and unstable manifolds are apparent from   figure \ref{fig_isosceles}.   Within the two-dimensional collision manifold $r=0$, the Lagragian restpoints are  saddles with one-dimensional stable and unstable manifolds.  In the three-dimensional energy manifold the starred restpoints pick up an extra stable dimension while the unstarred ones pick up  an extra unstable dimension.  Thus $W^u(L_\pm)$ and $W^s(L_\pm^*)$ are two-dimensional surfaces.  The Lagrangian homothetic connecting orbits represent transverse intersections of these surfaces.   Within the collision manifold,  $E_3^*$ is a repeller and $E_3$  is an attractor.
 Viewed within the three-dimesional energy manifold, $E_3$ and $E_3 ^*$  pick up an extra stable and unstable dimension, respectively.  Clearly the $E_3\rightarrow E_3^*$ connecting orbit, passing as it does through the interior of 
 the isosceles phase space,   is not a transverse intersection of stable and unstable manifolds.  This connecting orbit represents the Euler homothety orbit. 

A crucial fact for us is that the Eulerian restpoints have nonreal eigenvalues within the collision manifold for a large open set of masses, including those near equal masses.  This spiraling, together with the presence of restpoint connections in the collision manifold, causes the surfaces  $W^u(L_\pm)$ and $W^s(L_\pm^*)$ to wrap like scrolls around the $E_3\rightarrow E_3^*$ homothetic orbit as indicated in the figure.    To get a simpler picture of these scrolls,  we   use a piece of the plane $S=\{\nu=0\}$ as a cross-section to the flow near the homothetic orbit. The origin of $S$ represents the Euler homothetic orbit.  The   surfaces $W^u(L_\pm)$ and $W^s(L_\pm^*)$
 intersect $S$ in  spiraling curves.  More precisely (see \cite{M_InfinitelyClose}),   in polar coordinates $(\rho, \psi)$   near the origin,  each such curve is parametrized by the polar angle $\psi$
by setting  $\rho = f(\psi)$ where $f$ is some  strictly monotonic function: the curves wind monotonically around the origin.

\begin{figure}[h]
\scalebox{0.6}{\includegraphics{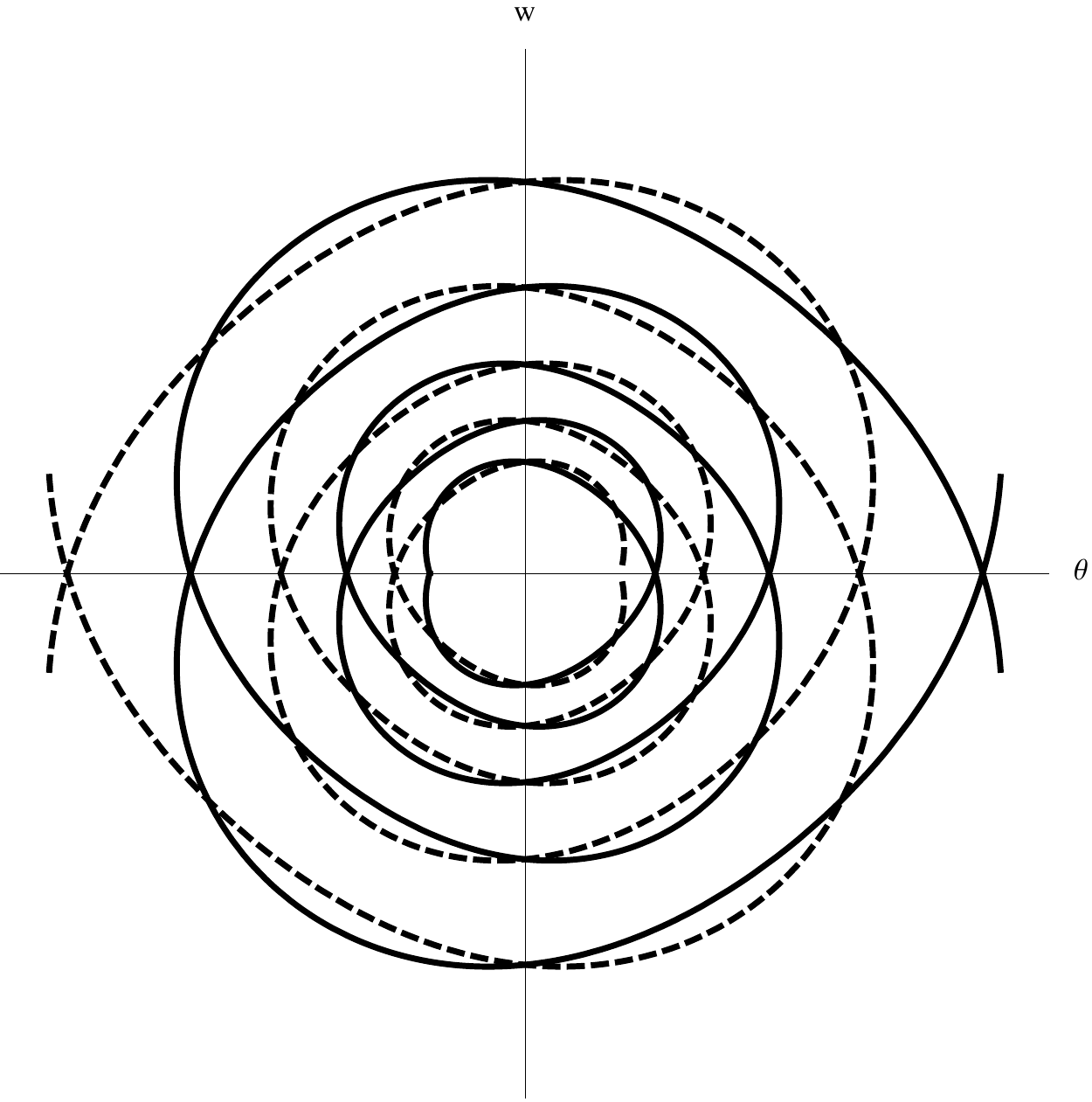}}
\caption{Spiraling stable and unstable manifolds. The stable manifolds $W^s(L_+ ^*), W^s(L_- ^*)$ spiral counterclockwise while the unstable manifolds $W^u(L_+), W^u(L_-)$ spiral clockwise.  The manifolds of $L_+$ are indicated with solid lines and those of $L_-$ are dashed. Each intersection point between stable and unstable  spiral represents  an isosceles solution doubly asymptotic to Lagrange triple collision. As we choose intersection points closer and closer to the origin,
the number of  syzygies  ($E_3$'s) suffered increases monotonically. }  
\label{fig_spirals}
\end{figure}

Figure~\ref{fig_spirals} shows the cross-section $S$ together with a schematic drawing of the four spiraling curves.  The curves representing $W^u(L_\pm)$ spiral in opposite directions from those representing $W^s(L_\pm^*)$.   The four curves are related by reflection symmetries.  Namely, let $\sigma_1(r,\theta,\nu,w) = (r,\theta,-\nu,-w)$ and $\sigma_2(r,\theta,\nu,w) = (r,-\theta,-\nu,w)$.  Then 
$\sigma_1(W^s(L_\pm^*))= W^u(L_\pm)$ and $\sigma_2(W^s(L_\pm^*))= W^u(L_\mp)$. 
It follows that there must be infinitely many intersection of each of the unstable manifolds with each of the stable ones.  Since the manifolds are real analytic, these intersections are either transverse or, at worst, finite-order tangencies.  Even a finite-order crossing is ``topologically transverse" and this is enough for the topological  construction given later on.  

We have found infinitely many topologically transverse intersections of Lagrangian ejection and collision orbits.  These orbits behave qualitatively as follows.  Beginning near one of the Lagrange ejection restpoints $L_\pm$ they follow a connecting orbit in the collision manifold to a neighborhood of the Eulerian ejection rest point $E_3$.  Then they follow the Eulerian homothetic orbit until they are close to the Eulerian collision restpoint $E_3^*$.  Finally, they follow another connecting orbit in the collision manifold to end at collision at one of the restpoints $L_\pm^*$

So far, this subsection has summarized previous work on the isosceles problem.  For the syzygy problem, we will need a few new observations.  The claim is that this type of orbit has a syzygy sequence of the form $3^n$ for some large value of $n$ and, moreover, any sufficiently large $n$ can be achieved in this way.  To see this,  first note that the connecting orbits $L_\pm\rightarrow E_3$ and $E_3^*\rightarrow L_\pm^*$ in the collision manifold go directly from the equilateral to collinear shapes, staying far from double collision and far from any syzygy other than type $3$.  

To study syzygies of type $3$ recall that the collinear shape is represented in our $(r, \theta, w, \nu)$ coordinates by the plane $\theta =  0$, which in the plane $S$ of figure \ref{fig_spirals}
is represented by the $w$-axis. The  plane $\theta =0$ contains the Eulerian restpoints and  Eulerian homothetic orbit.  Due to the spiraling at the restpoints, it is clear that our orbits will have a large number of type $3$ syzygies.  It is also clear that the number is odd for connections $L_+\rightarrow L_-$ and $L_-\rightarrow L_+$.  It remains to show than any sufficiently large number of syzygies can be achieved.

One of the differential equations for the blown-up isosceles problem is $\theta' = w$.  Thus, except for the Euler homothetic orbit, orbits cross the syzygy plane transversely.  In figure~\ref{fig_isosceles} orbits in the front half of the energy manifold ($w>0$) are crossing from right to left while those at the back ($w<0$) cross from left to right.    It follows that the number of syzygies which occur on any orbit segment is a continuous function of the endpoints of the segment, as long as the endpoints are not on the syzygy plane.

Now consider one of the two stable spirals in figure~\ref{fig_spirals}.  Each point of the spiral determines a forward orbit segment ending, say, in some convenient cross-section near the corresponding Lagrange restpoint.   It follows that as we vary the point on the spiral,  the total number of syzygies experienced in forward time can change only when the initial point crosses the line $\theta = 0$.  Since the spiraling is monotonic, the plane $\theta = 0$ divides the spiral into curve segments in each of which the total number of forward syzygies is constant.   Moreover, points near $\theta = 0$ immediately cross it once or have just crossed it.  Hence for nearby points of the spiral on opposite sides of $\theta = 0$, the total number of syzygies differs by one.  Since the number of syzygies tends to infinity as the spiral converge to the center, every sufficiently large number of forward syzygies is attained.  There is a similar story for counting backward syzygies of points on the unstable spirals.

Consider a segment $\gamma_n$ of $W_s(L_+^*)$ with endpoints in $\theta=0$ whose interior points have exactly $n$ forward time syzygies.  Then $\sigma_1(\gamma_n)$ is a segment of $W^u(L_+)$ whose interior points have $n$ syzygies in backward times.  These spirals intersect at a point of the line $w=0$ (the zero velocity curve).  This point represents one of our ejection-collision orbits having exactly $2n$ syzygies of type $3$.  Similarly there are $L_-^*\rightarrow L_-$ ejection collision orbits with exactly $2n$ syzygies.

On the other hand, for the same segment $\gamma_n$ of $W_s(L_+^*)$ the other reflection $\sigma_2(\gamma_n)$ is a segment of $W^u(L_-)$.  These intersect at their endpoints, that is, along $\theta=0$.  These intersection points represent $L_+^*\rightarrow L_-$ ejection collision orbits.   One endpoint or the other  has exactly $n$ syzygies in both forward and backward time including the initial syzygy (if the segment lies in $\theta\le 0$ it's the endpoint with $w<0$; otherwise the one with $w>0$).
For this endpoint, the total number of forward and back ward syzygies will be exactly $2n-1$.  Thus any sufficiently large number of syzygies can be achieved.

\subsection{Symbolic dynamics in the planar problem}
Consider the planar three-body problem with equal masses.  For each fixed angular momentum $\mu$, we get a rotation-reduced Hamiltonian system of three degrees of freedom. 
When the energy is further fixed at the value $h$ we get a flow on a 5-manifold denoted $M(h,\mu)$.
 If the angular momentum is zero then three separate isosceles problems appear as subsystems.  
 For each of these we can construct ejection-collision orbits as in the last subsection.  So we can realize syzygy sequences of types $1^n, 2^n$ and $3^n$ for all sufficiently large $n$,  say $n\ge N$.  In this subsection we will show how to concatenate these blocks to get bi-infinite syzygy sequences of the form $\ldots \epsilon_{-1}^{j_{-1}}\epsilon_0^{j_0}\epsilon_1^{j_1}\ldots$ where $\epsilon_i\in\{1,2,3\}$ and $n_i\ge N$.  This requires perturbing to nonzero angular momentum and for this we need to choose an upper bound $N^{\prime}$ for $n_i$ as in the statement of theorem~\ref{th_2}.   In addition we will see that nonequal masses can be handled and that every periodic sequence of this form is realized by at least one periodic orbit.  The rest of this section is a summary of the results of \cite{M_chaotic}  which we refer for more details.  See also \cite{M_symbolic, M_heteroclinic} for some similar arguments.
Our aim   is to convey the spirit of the proof and not all the details, and to point out where a bit of new work is needed to keep track of stutter block lengths.

The McGehee blow-up for the planar problem  has a lot in common with the blow-up for the isosceles problem.  The energy manifolds $M(h,0)$
 are five-dimensional and  all share  the four-dimensional  collision manifold $r=0$ as a common boundary.   There are five central configurations which lead to a total of ten restpoints in the collision manifolds, namely five ejection restpoints $E_1,E_2,E_3,L_+,L_-$ and the corresponding collision restpoints with starred notation.  The stable and unstable manifolds $W_s(L_+^*)$ and $W_u(L_+)$ now have dimension three and the Lagrange homothetic solutions are still transverse intersections of these.  Moreover, all of their infinitely many other topologically transverse intersections viewed within the isosceles submanifolds remain  topologically transverse when viewed from  the larger planar energy manifolds.  There are two extra dimensions to be understood in the planar problem.  The cross-section $S=\{\nu=0\}$ (figure ~\ref{fig_spirals}) to the Euler homothetic orbit is now four-dimensional and instead of spiraling curves we have spiraling two-dimensional surfaces.  It is shown in \cite{M_chaotic} that these surfaces intersect topologically transversely in $S$,

The realization of arbitrary concatenations of stutter blocks is accomplished by the methods of symbolic dynamics analogous to those used for the Smale horseshoe map. (See \cite{Guckenheimer_Holmes}.)  One version of the horseshoe is built by taking two rectangles each of which is stretched in one direction and contracted in the other and then mapped across both of the original rectangles.  Label the two rectangles by symbols $0,1$.  One can   prescribe an arbitrary bi-infinite sequence of $0$'s and $1$'s and prove that the sequence is realized by a unique orbit.  The sequence itself  is   the {\em itinerary} of the orbit: the list of rectangles visited in order of visit.  The uniqueness comes from uniform hyperbolic stretching.  In the topological approach used here, we will not be able to guarantee uniqueness.

\begin{figure}[h]
\scalebox{0.6}{\includegraphics{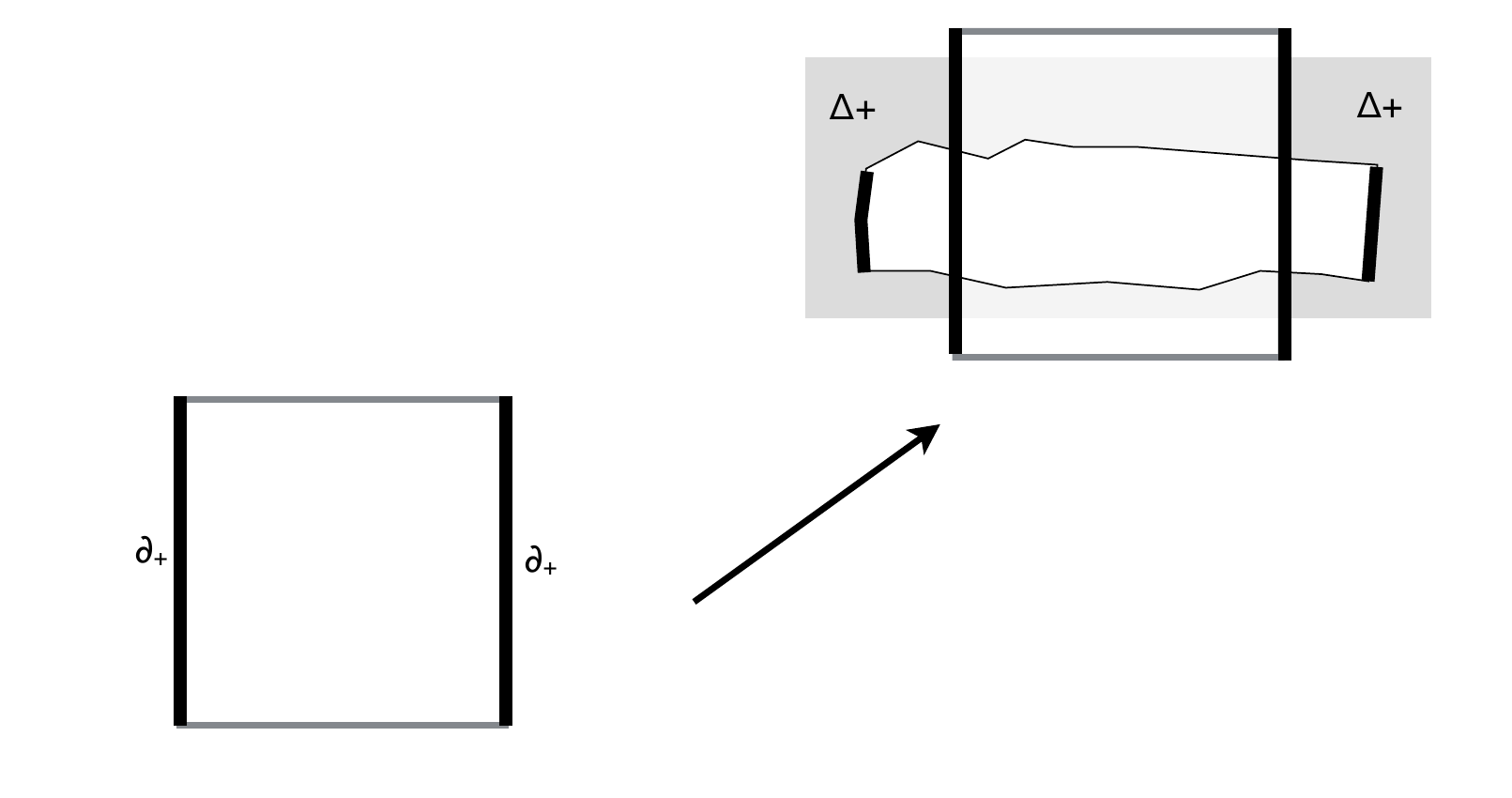}}
\caption{Two dimensional windows being correctly aligned by a Poincar\'e map.  The first window maps into an auxiliary window set up near the target window. }  
\label{fig_windowpic}
\end{figure}

The topological approach based on ``windows" was pioneered by Easton \cite{Easton} and further   developed and used in \cite{M_chaotic}  for the planar three-body problem. We follow the discussion and notation of \cite{M_chaotic},  p.56-60, closely.   Instead of the rectangles in the horseshoe map, we have four-dimensional boxes called windows homeomorphic to the product $D_+ ^2\times D_{-}^2$ of two two-dimensional disks. The analogous splitting for a two-dimensional window is shown in figure~\ref{fig_windowpic}.  In that figure, horizontal line segments of one window are being stretched horizontally across the next box.  For our four-dimensional window $w\simeq D_+ ^2\times D_{-}^2$, the first disk represents two directions which will be stretched in forward time and will be called positive.  The second disk represents two directions which will be  compressed in forward  time, or alternatively, stretched  in backward time. These directions will be called negative.  There is a splitting of the boundary of  $D_+ ^2\times D_{-}^2$  into two parts $\partial^+ = \partial D_+ ^2\times D_{-}^2$ and $\partial^-=D_+ ^2\times \partial D_{-} ^2$, each homeomorphic to a solid torus.   A positive disk has its boundary in $\partial_+$ and represents a generator of the relative homology group $H_2(w,\partial^+)\simeq \mathbb{Z}$.

  Such a window can be constructed near each of the transverse intersections of the surfaces  $W^s(L_\pm^*)$ and $W^u(L_\pm)$ in the four-dimensional cross-section $S$.   In carrying out the  perturbation to nonzero angular momentum later on, we will only be able to work with finitely many windows.    This is where the choice of the number $N^{\prime}$ in the statement of theorem~\ref{th_2} comes in.   We already have $N$ as lower bound on the number of syzygies which are realized by the ejection collision orbits.  Given any $N^{\prime}\ge N$ we can restrict attention to the finitely many topologically transverse ejection-collision orbits which realize stutter blocks size in the range $N\le n\le N^{\prime}$.  Then our perturbation will be constructed based on this finite collection of windows.

\begin{figure}[h]
\scalebox{0.6}{\includegraphics{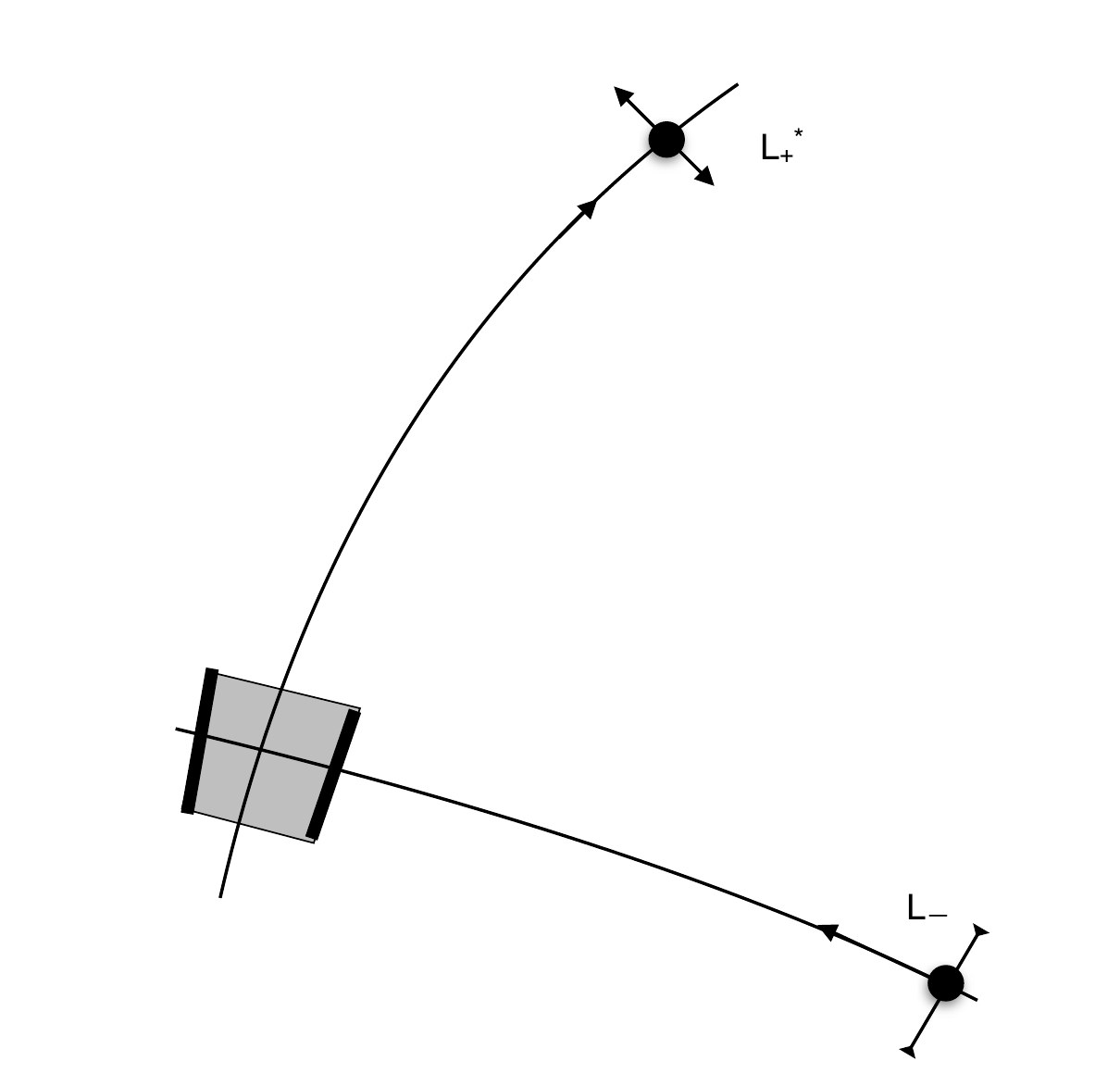}}
\caption{Two-dimensional representation of how to select a window.  As in the four-dimensional construction,  the positive boundary $\partial_+$ (bold edges) is linked with the stable manifold of the restpoint and will be stretched in forward time.  This shows the situation before perturbing to nonzero angular momentum.  We actually use nearby windows in the nonzero angular momentum manifolds.}
\label{fig_saddlewindow}
\end{figure}

The windows are chosen so that the positive directions are aligned with the unstable manifold involved in the intersection while the negative directions are aligned with the stable manifold (thus a positive disk crosses the stable manifold transversely in some sense and, when followed forward near the restpoint, will get stretched; the positive boundary $\partial_+$ is linked with the stable manifold).  For definiteness,  consider a point of intersection of $W^u(L_-)$ and $W^s(L_+^*)$ which realizes the stutter block $3^{n}$ for some odd number $n$.  Using topological transversality, one can choose a $C^0$ local coordinate system which makes these manifolds coordinate planes and then choose the window to be a product of small discs in these planes.  See figure \ref{fig_saddlewindow}. If the window is sufficiently small, the orbits starting there will behave qualitatively like the ejection-collision orbit at its center.  They will enter a small neighborhood of $L_+^*$ in forward time and a small neighborhood of $L_+$ in backward time.  In between they will avoid collisions, realize the stutter block $3^n$, and stay close to the isosceles manifold for $3$.    Setting up one such  window near 
each of a finite number of intersection points whose corresponding solutions have  $n$ stutters with $n\in[N,N_1]$, 
we need to show that when we perturb to small nonzero angular momentum, each of these windows is stretched across each of the others by the perturbed flow.  A topological analogue of the uniform stretching in the Smale horseshoe is provided by a homological definition of  ``correct alignment".  The version used in \cite{M_chaotic} differs from that in Easton's work and will be described briefly now.

We can't expect the Poincar\'e mappings of the three-body flow to carry one window exactly onto another.  Rather, as in the horseshoe, a window is stretched through the other and overlaps it.  For this reason it is convenient to set up some larger auxiliary windows which capture the homology of the original as shown schematically in figure~\ref{fig_windowpic}.  To each window $w$ we associate an auxiliary window $W^+$ with a thickened boundary set $\Delta^+= W^+\setminus w$ such that there is a retraction map $r_+:(W^+,\Delta^+)\into (w,\partial^+)$ fixing $w$ and inducing an isomorphism on relative homology groups $H_2(W^+,\Delta^+)\simeq H_2(w,\partial^+)$.  
Similarly, each window will have an associated auxiliary pair $(W^-,\Delta^-)$.  Then two windows $w_0, w_1$ are said to be {\em correctly aligned} if there is a flow-induced Poincar\'e map taking $(w_0,\partial_0^+) \into (W_1^+,\Delta_1^+)$ and inducing an isomorphism on the second relative homology groups.  We also require the inverse Poincar\'e map to take  $(w_1,\partial_1^-) \into (W_0^-,\Delta_0^-)$ and to induce an isomorphism on homology.  So if two windows $w_0,w_1$ are correctly aligned there is a kind of homological stretching instead of the usual hyperbolic stretching.   


The next step is to  show that given any bi-infinite sequence of correctly aligned windows, there is a nonempty compact set of orbits which passes through their interiors using the given Poincar\'e maps.  To get the idea of the proof,   recall the proof of the analogous statement for the Smale horseshoe.  Consider any bi-infinite itinerary specifying which of the two rectangles the orbit should hit at each iteration.  The set of points in the initial rectangle which map into the next one is a negative subrectangle of the initial rectangle.   In fact any finite forward time itinerary is realized by such a subrectangle.  Similarly, any finite backward time itinerary is realized by a positive subrectangle.  By intersecting these rectangles we get a nonempty compact set realizing any finite itinerary.  Since an intersection of nested, nonempty compact sets is nonempty and compact, we can also realize a bi-infinite itinerary. From the topological perspective,   we replace the idea of negative subrectangle with ``compact set" which intersects every positive segment' or better ``compact set which intersects the support of every chain which is nontrivial in $H_1(w,\partial^+)$".  

  Returning to our four-dimension case, define a {\em positive chain} to be a relative singular two-chain whose homology class in $H_2(w,\partial^+)$ is nonzero and similarly for {\em negative chains}.   Using the definition of correct alignment, one can show that the set of points in the initial box which realize any given finite forward time itinerary is a compact set which intersects the support of any  positive chain.  
  Similarly the points realizing a chosen finite backward itinerary form a compact set  intersecting the support of every negative chain.   Using some algebraic topology one can show that any two such  sets of this type have  nonempty intersection and the proof for bi-infinite sequences is completed as for the horseshoe.

The proof that a periodic sequence is realized by at least one reduced periodic orbit is similar.  Here we have a composition of Poincar\'e maps giving a map $\tilde \phi$ from a subset of the initial box to itself.  Consider $F_1=\{(x_1,x_2)\in D_+ ^2\times D_{-}^2: \tilde\phi_1(x_1,x_2)=x_1\}$, that is, the set of points whose first two coordinates (out of four) are fixed.   Similarly let $F_2=\{(x_1,x_2)\in D_+ ^2\times D_{-}^2: \tilde\phi_2(x_1,x_2)=x_2\}$.  Then one can show  that $F_1$ intersects the support of every positive chain and $F_2$ intersects the support of every negative chain, and it follows as above that $F_1\cap F_2\ne\emptyset$, that is, there is at least one fixed point.

\begin{figure}
\scalebox{0.8}{\includegraphics{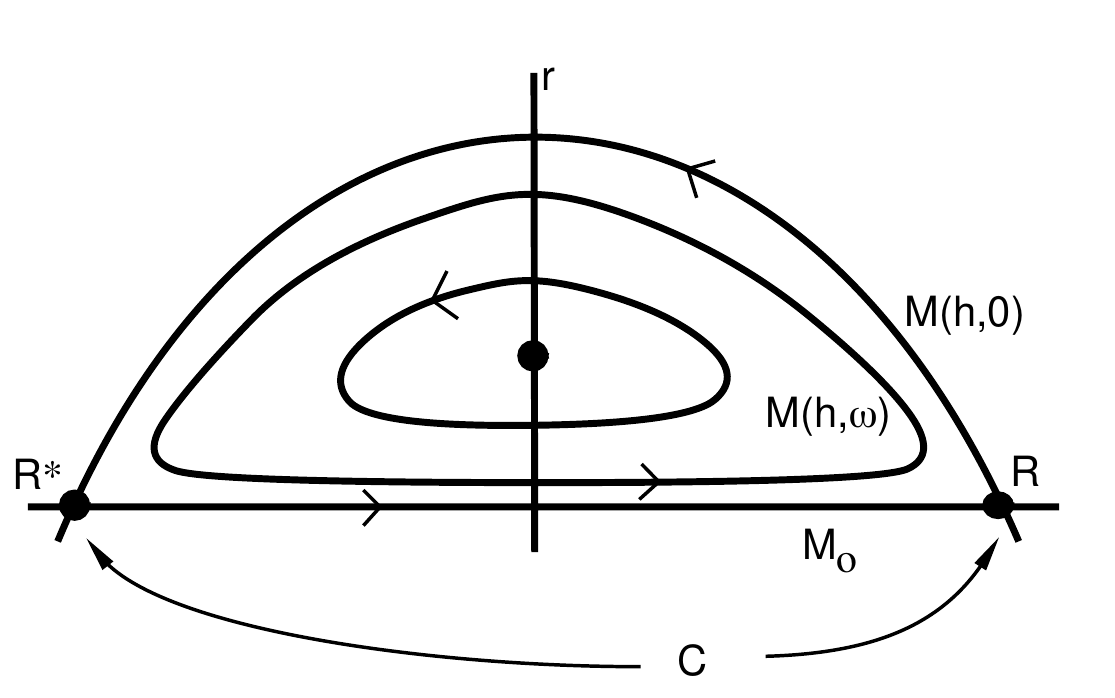}}
\caption{A cycle of restpoint connections.  $R$ represents one of the restpoints $L_+, L_-$.  The connecting orbit $R\into R^*$ is one of the the Lagrange homothetic orbits while the connection $R^* \into R$ in \{r=0\} represents the limit as angular momentum tends to zero of the spinning of the Lagrange homographic solutions near triple collision.  The figure also shows schematically how the nonzero angular momentum manifolds converge to a union of the zero angular momentum manifold and another five manifold $M_0$ at $r=0$. }  
\label{fig_loop}
\end{figure}

It  remains to explain how we can get each of the  windows defined near our ejection-collision orbits to be correctly aligned with itself and all of the other windows when the angular momentum is nonzero and sufficiently small.   We know that each   window  approaches one of the hyperbolic collision restpoints $L_\pm^*$ in forward time and one of the ejection restpoints $L_\pm$ in backward time.   We need to find a way to go from a neighborhood of $L_\pm^*$ to a neighborhood of $L_\pm$.   It turns out that there are transverse connecting orbits $L_+^*\into L_+$ and $L_-^*\into L_-$ contained in the boundary flow at $r=0$.  We have already discussed  the four-dimensional collision manifold at $r=0$ which serves as a boundary to the $\mu=0$, fixed energy manifolds.  But there is another five-dimensional flow at $r=0$ representing the limits as $\mu\into 0$ of nonzero angular momentum orbits ($M_0$ in figure~\ref{fig_loop}).  In particular, consider the highly eccentric Lagrange homographic solutions.  
The equilateral triangle formed by the three bodies expands and contracts to near zero.  Instead of colliding, the triangle quickly spins by $360^\circ$ and the triangle expands again.  In the limits as $\mu\into 0$ these solutions converge to  cycles of restpoint connections.  We get the Lagrange homothetic orbits $L_\pm \into L_\pm^*$ and also connections going in the other direction  $L_\pm^* \into L_\pm$.  The latter represent the limit of the spinning behavior.  These new connections are also transverse and we set up another window along each of them.

\begin{figure}[h]
\scalebox{0.6}{\includegraphics{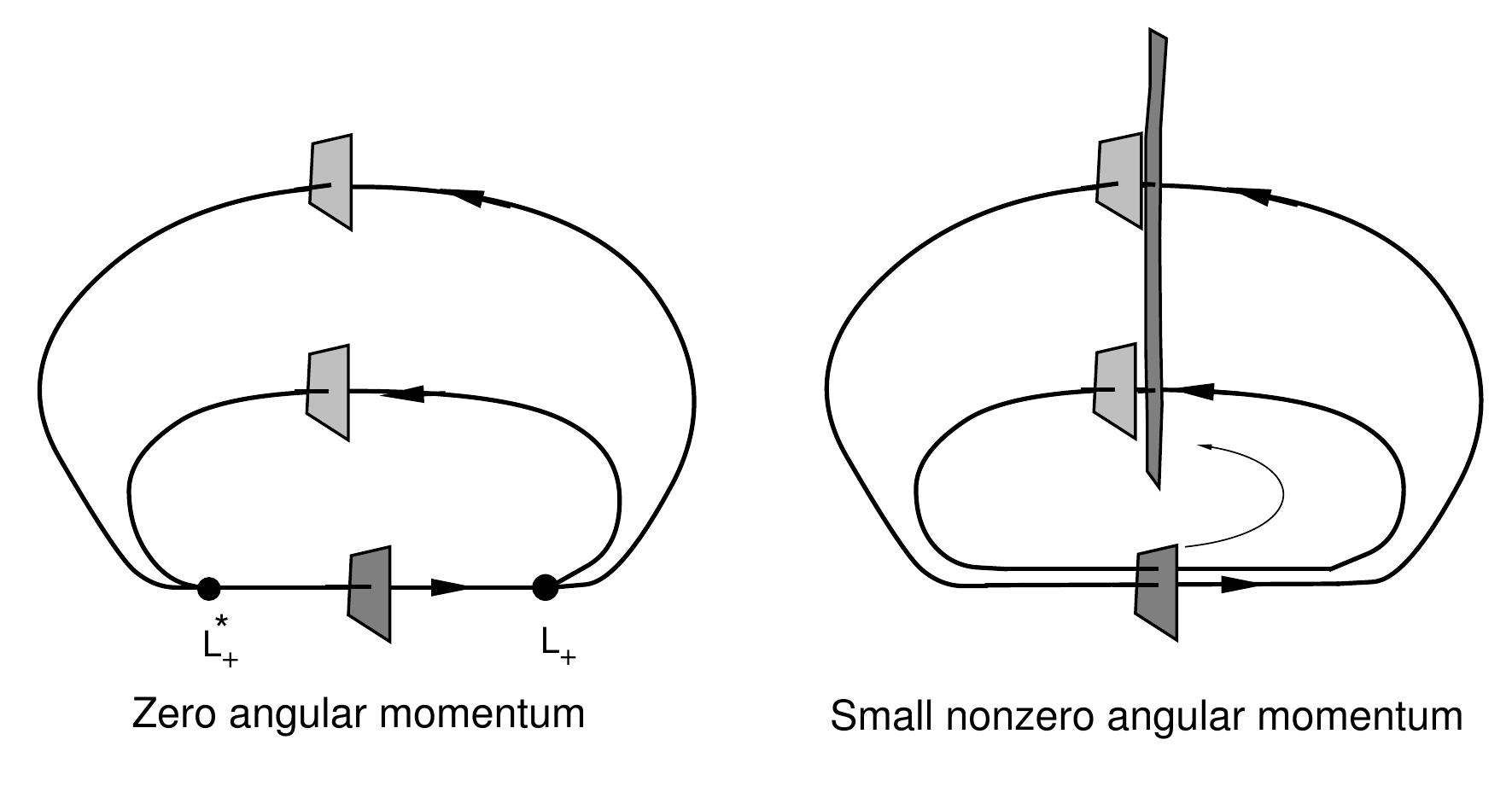}}
\caption{Perturbing to nonzero angular momentum removes the restpoints and allows the windows to align with one another.}  
\label{fig_perturb}
\end{figure}

For $\mu=0$ the windows near our ejection-collision orbits approach the Lagrangian restpoints and are stretched out by the hyperbolic dynamics there.  Similarly the new windows constructed along the connecting orbits in $\{r=0\}$ approach these same restpoints and are also stretched out.  But the four-dimensional collision manifold $C$ stands between them
within $M_0 \cup M(h,0)$  and prevents them from aligning under the flow.  It turns out that  when we perturb   to nonzero angular momentum,  the restpoints disappear, 
as does $C$ itself, allowing each incoming window to flow across the neighborhood of the former restpoints to meet each outgoing window.    This is shown schematically in figure~\ref{fig_perturb}.  Thus there is a chain of alignments connecting each of our ejection-collision windows with each of the others via the new window near $r=0$.  Since the qualitative behavior represented by the new window is just the spinning of a small  equilateral triangle, no new syzygies or collisions are introduced.

\section{Open Questions}

1.  Does theorem 1  hold for zero angular momentum?    
This   was our original question,   suggested by Wu-Yi Hsiang. 
Solutions of Theorem 2  converge to   triple ejection-collision  solutions when they are continued to angular momentum zero.
Another   idea will be needed to yield such realizing solutions.   

2.  By how much can we increase angular momentum before our  realizing classes  disappear? 
Eventually
all of them   disappear for topological reasons except those realizing the ``tight binary" 
reduced  free homotopy classes, meaning the syzygy sequences $12$, $13$, and $23$ repeated periodically.  
Topologically  speaking,  all free homotopy classes are   possible below the ``Euler value" in the equal mass case.  (See \cite{M_contemporary}.)  
  Do  realizations  survive right up to this Euler value?
 
3. Does near-isosceles imply isosceles?  Suppose that the angular momentum is zero,  the energy is negative
and the masses are all equal.  Is  there  a positive $\epsilon$ such that any solution which spends its
whole existence within $\epsilon$ of the isosceles submanifold must be
an isosceles solution?  Here ``within $\epsilon$" means as measured in the shape sphere,
computing the distance of the projected solution to the isosceles great circles.   This question is motivated by contemplating the likely fate of the solutions described
in this paper as the angular momentum tends to zero, at fixed negative energy. 

{\bf Remark.} The analogous
question with  the collinear submanifold replacing the isosceles submanifold has the answer ``no".
The Schubart collinear orbit is   KAM stable   with respect
to noncollinear perturbations.  Near to it we have persistent torii filled with noncollinear   solutions.
 
4.  Can we get rid of all the stutters in a solution realizing a given reduced free homotopy class? 
The realizing solutions of Theorem 2 have lots of stutters.  For the equal mass, under the strong force assumption,  and in the case of  angular momentum  zero, 
we have proved (\cite{Mont_pants}) that   realizing solutions have no stutters.  (This no-stuttering  follows from  
the fact that  the realizers in this case  are variational minimizers  and variational minimizers cannot have stutters.)
More generally, and less precisely, how far can  periodic solutions   be pushed
away from triple collision,   away from the Hill boundary, into the ``interior" of the three-body problem?

5.  We can use the homotopy method to  follow  any realizing solution from the equal mass,  strong force,  angular momentum case as  discussed in the previous question,  changing   the potential, the energy, and the  angular momentum continuously along some short arc.
Do any of these persist all the way to the case of the  Newtonian potential, and in so doing   connect up with
the realizing solutions described here by Theorem 2?  

6.   Which  reduced free homotopy classes admit dynamically KAM stable realizing solutions?
All   realizing solutions obtained by Theorem 2 here    are unstable.
We do not    know many  such ``stable  classes". An initial perusal yields only the class $123123$ of the figure eight,   the classes $12$, 
$13$ and $23$  as realized by  planetary configurations (depending on which is the dominant mass) and earth-moon-sun type configurations,
and the class $1232$  realized by  Broucke-Henon  solutions.

\section{acknowledgements}  We would like to thank Wu Yi Hsiang for posing a version of the realization 
problem back in 1997.  Thanks to Carles Simo for challenging us to look for a dynamical
mechanism behind realization.  R.Moeckel and R.Montgomery  gratefully acknowledge  NSF grants DMS-1208908 and DMS-1305844, respectively,   for essential support.


\

\end{document}